\documentclass[12pt]{article}
\usepackage{amsmath,amssymb,amsbsy}
\usepackage{graphicx}
\usepackage{diagbox,multirow}
\usepackage{tikz}
\usepackage{tkz-graph}
\usetikzlibrary{shapes.geometric}

\usepackage{tkz-euclide,amsmath} 
\usepackage{soul}
\usepackage{url}
\usepackage{subcaption}
\usepackage{caption}
\usepackage{thm-restate}

\newtheorem{theorem}{Theorem}
\newtheorem{lemma}[theorem]{Lemma}
\newtheorem{corollary}[theorem]{Corollary}
\newtheorem{proposition}[theorem]{Proposition}

\newtheorem{observation}[theorem]{Observation}

\newtheorem{definition}[theorem]{Definition}

\newtheorem{conjecture}[theorem]{Conjecture}
\newtheorem{open}[theorem]{Open Problem}

\newtheorem{remark}[theorem]{Remark}

\newcommand{\qed}{\hfill \rule{.1in}{.1in}}

\def\mod{\mathop{\rm mod}\nolimits}
\def\zet{\mathop\mathbb{Z}\nolimits}

\makeatletter
\def\imod#1{\allowbreak\mkern10mu({\operator@font mod}\,\,#1)}

\title{\bf Open XOR-magic odd graphs and closed XOR-magic even graphs}

\begin{document}

\author{{{Sylwia Cichacz$^{1}$\thanks{The work of the author was supported by  the AGH University of Krakow under grant no. 16.16.420.054, funded by the Polish Ministry of Science and Higher Education.}, Hubert Grochowski$^2$\thanks{Hubert Grochowski’s research was funded by the Warsaw University of
Technology within the Excellence Initiative: Research University (IDUB)
programme.}, Rita Zuazua$^{3}$}}
	\footnote{e-mails: cichacz@agh.edu.pl, hubert.grochowski.dokt@pw.edu.pl,\\ ritazuazua@ciencias.unam.mx}\\
\normalsize $^1$AGH University of Krakow, Poland \\
\normalsize $^2$Warsaw University of Technology, Poland\\
\normalsize $^3$National Autonomous University of Mexico, Mexico
}

\maketitle
\begin{abstract}
XOR-magic graph labelings form a special subclass of group distance magic labelings. A simple connected graph of order $2^n$
 is called an \textit{open} (respectively, \textit{closed}) \textit{XOR-magic graph of power $n$} if its vertices can be labeled bijectively with vectors from $(\zet_2)^n$  such that the sum (over $(\zet_2)^n$) of labels in each open (respectively, closed) neighborhood of every vertex is equal to the zero vector.

In \cite{Batal}, Batal asked whether there exists any odd-regular open XOR-magic graph or any even-regular closed XOR-magic graph. In this paper, with partial help of MILP solver, we answer this question in the affirmative. More precisely, we prove that for every integer $n>3$, there exists an odd-regular open  XOR-magic graph of power $n$ and an even-regular closed XOR-magic graph of power 
$n$.  We also show some applications of the spectra of graphs for an open XOR-magic labeling.
\end{abstract}


\section{Introduction}

Let $G=(V,E)$ denote a graph with vertex set $V$ and edge set $E$. For {$x \in V$} we denote the set of all neighbors of $x$ by {$N(x)$} and $N[x]=N(x)\cup\{x\}$. The set $N(x)$ we call \textit{open neighborhood set} of $x$, whereas $N[x]$ we call \textit{closed neighborhood set} of $x$. If there exists an integer $k$ such that $|N(x)|=k$ for every $x \in V$, then we say that $G$ is {\em regular with valency} $k$, {or $G$ is a $k$-regular graph}. Furthermore, we say that $G$ is {\em odd} ({\em even}) {\em regular}, if it is regular with valency $k$, where $k$ is odd (even, respectively).

Let $(\Gamma,+)$ be a finite Abelian group.  A $\Gamma$\emph{-distance magic labeling} of a graph $G = (V, E)$ with $|V| = n$ is a bijection $\ell$ from $V$ to an Abelian group $\Gamma$ of order $n$, for which there exists $\mu \in \Gamma$, 
such that the weight $w(x) =\sum_{y\in N(x)}\ell(y)$ of every vertex $x \in V$ is equal to $\mu$ \cite{Fro}. In this case, the element $\mu$ is called the \emph{magic constant of} $G$ and the graph $G$ is called $\Gamma$-distance magic graph. Note that in a $\Gamma$-distance magic labeling, the weight of each vertex $v$ is computed over its open neighborhood $N(v)$. Analogously, in \cite{ACFSQ}, the notion of a $\Gamma$-closed distance magic labeling is defined, where the weight is taken over the closed neighborhood $N[v]$. Formally, a $\Gamma$\emph{-closed distance magic labeling} is a bijection $\ell'\colon V \to \Gamma$, such that the weight $w'(v) =\sum_{u\in N[v]}\ell(u)$ of every vertex $v \in V$ is equal to the same element $\mu'\in \Gamma$. 

Independently, Siehler defined a closed XOR-magic graph of power $n$. Formally speaking, a simple connected graph of order $2^n$ is defined as a \textit{closed XOR-magic graph of power} $n$  if its vertices can be labeled with vectors from $(\zet_2)^n$ in a one-a
to-one manner, such that the sum of labels in each closed neighborhood set
of vertices equals zero \cite{Siehler}. Thus, a closed XOR-magic graph of power $n$ is a connected  $(\zet_2)^n$-closed distance magic graph with the magic constant $0$. A labeling ensuring the fulfillment of the condition by definition will be called a \textit{closed XOR-magic labeling}.

This definition was recently adopted by Batal for an open neighborhood of vertices, an open XOR-magic graph of power $n$ is a connected $(\zet_2^n)$-distance magic graph with the magic constant $0$ \cite{Batal}. Such labeling we call an \emph{open XOR-magic labeling}.  Batal proved the
existence of $k$-regular closed XOR-magic graphs of power $n$ for all $n\geq 2$ and for all
$k\in \{3,5,7,\dots, 2^n\text{-}5\}\cup \{2^n\,\text{-}\,1\}$. He  showed that there is no
($2^n\,\text{-}\,3$)-regular open XOR-magic graph of power $n$. 

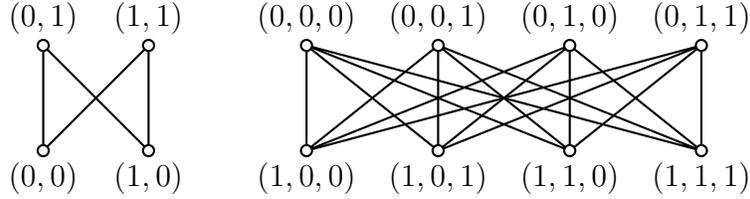
\begin{figure}[ht!]
\begin{center}
\begin{tikzpicture}[scale=0.7,style=thick,x=1cm,y=1cm]
\def\vr{3pt} 

\path (0,0) coordinate (a1);
\path (2,0) coordinate (b1);
\path (2,2) coordinate (c1);
\path (0,2) coordinate (d1);

\draw (c1) -- (b1) -- (d1) -- (a1) -- (c1);

\draw (a1) [fill=white] circle (\vr);
\draw (b1) [fill=white] circle (\vr);
\draw (c1) [fill=white] circle (\vr);
\draw (d1) [fill=white] circle (\vr);

\draw[anchor = north] (a1) node {$(0,0)$};
\draw[anchor = north] (b1) node {$(1,0)$};
\draw[anchor = south] (c1) node {$(1,1)$};
\draw[anchor = south] (d1) node {$(0,1)$};


\path (5,0) coordinate (a);
\path (7.5,0) coordinate (b);
\path (10,0) coordinate (c);
\path (12.5,0) coordinate (g);
\path (5,2) coordinate (d);
\path (7.5,2) coordinate (e);
\path (10,2) coordinate (f);
\path (12.5,2) coordinate (h);

\draw (a) -- (d) -- (c) -- (e) -- (a) -- (f) -- (b) -- (e);
\draw (b) -- (d) -- (g) -- (e);
\draw (c) -- (f) -- (g) -- (h) -- (c);
\draw (a) -- (h) -- (b);

\draw (a) [fill=white] circle (\vr);
\draw (b) [fill=white] circle (\vr);
\draw (c) [fill=white] circle (\vr);
\draw (d) [fill=white] circle (\vr);
\draw (e) [fill=white] circle (\vr);
\draw (f) [fill=white] circle (\vr);
\draw (g) [fill=white] circle (\vr);
\draw (h) [fill=white] circle (\vr);

\draw[anchor = north] (a) node {$(1,0,0)$};
\draw[anchor = north] (b) node {$(1,0,1)$};
\draw[anchor = north] (c) node {$(1,1,0)$};
\draw[anchor = north] (g) node {$(1,1,1)$};
\draw[anchor = south] (d) node {$(0,0,0)$};
\draw[anchor = south] (e) node {$(0,0,1)$};
\draw[anchor = south] (f) node {$(0,1,0)$};
\draw[anchor = south] (h) node {$(0,1,1)$};

\end{tikzpicture}
\end{center}
\caption{A $(\zet_2)^2$-distance magic labeling of $K_{2,2}$ and a $(\zet_2)^3$-distance magic labeling of $K_{4,4}$ }
\label{sporadic}
\end{figure}\label{DMXOR}

Observe that in Figure~\ref{DMXOR}, the $(\mathbb{Z}_2)^2$-distance magic labeling of $K_{2,2}$ is not an open XOR-magic labeling, since its magic constant $\mu=(1,0)$ is different from $(0,0)$. In contrast, the $(\mathbb{Z}_2)^3$-distance magic labeling of $K_{4,4}$ is an open XOR-magic labeling.

\vspace{.3cm}

Batal showed some connection for open XOR-magic and closed XOR-magic graphs.
 \begin{observation}[\cite{Batal}]\label{complement} Let $G$ be a connected   graph on $m=2^n$ vertices such that $\overline{G}$ is also connected. Then $G$ is an open XOR-magic graph if and only if $\overline{G}$ is a closed XOR-magic graph.
 \end{observation}

While there are many papers on even regular $\Gamma$-distance magic graphs see  \cite{CicFro,DeGeZe} for more references, there is not much known for odd regular:
\begin{theorem}[\cite{CicFro}]\label{gr:odd} Let $G$ be an $r$-regular graph on $n$
 	vertices, where $r$ is odd. If $\Gamma$ is an Abelian group of order $n$ with $|I(\Gamma)|=1$, then $G$ is not $\Gamma$-distance magic
 \end{theorem}
 
 \begin{observation}[\cite{CicFro}]\label{obs:odd} Let $G$ be an $r$-regular  graph on $n\equiv 2 \pmod4$ vertices, where $r$ is odd. Then $G$ is not $\Gamma$-distance magic for any Abelian group $\Gamma$ of order $n$.
 \end{observation}


However, there are some known examples of odd regular $\Gamma$-distance magic graphs \cite{Cic,CicMik}. Only three of them are connected, and   there were no known ones for $\Gamma=(\mathbb{Z}_2)^m$. Therefore, Batal posted the following open problem:
\begin{open}\label{XORopenproblem}
Is there any even regular closed XOR-magic or odd regular open
XOR-magic graph?
\end{open}
In this paper, we answer this question affirmatively in {Section~\ref{odd_even}}. Namely, we prove the following theorem.
\begin{restatable}{theorem}{CoPthm}
  \label{th:C1P}
For any integer $n>3$ there exists an even closed XOR-magic graph of power $n$, and an odd open XOR-magic graph of power $n$.
\end{restatable}
Moreover, in Section~\ref{algebraic} we show some applications of the spectra of graphs for open XOR-magic labeling.

\section{Algebraic approach and its application to Cayley graphs}\label{algebraic}

In this section, we focus on the algebraic approach to the problem of the nonexistence of open XOR-magic labelings for graphs. Then, we will apply our approach to Cayley graphs.
Let $X$ be  a  (finite or infinite) group
and  $S$  be  a  unit-free  symmetric subset of $S$. The \textit{Cayley graph} $\mathrm{Cay}(X,S)$  is a graph with vertices labeled (bijectivelly) by the elements of the group $X$, and the adjacency relation is defined by the subset $S$  of the group $X$.  In particular, \textit{circulant graphs} are Cayley graphs defined over cyclic groups -- the simplest and most fundamental family of groups. Formally, let $m \geq 4$, $D = \{1,2,\ldots \lfloor{m/2}\rfloor\}$ and let $S = \{s_1,s_2,\ldots,s_r\} \subseteq D$. Then, the \emph{circulant graph} $G = C_m(S)$ has $m$ vertices $x_0, x_1, \ldots x_{m-1}$ in which two distinct vertices $x_i$ and $x_j$ are adjacent if and only if $|i - j| \in \{s_i, n - s_i : s_i \in S\}$. 

\vspace{.3cm}
In the following, for each $a, b \in \mathbb{Z}$ such that $a < b$, by $[a,b]$ we denote the set $\{a,a+1,\ldots,b\}$.

By the definition of open (closed) XOR-magic graph, we will concentrate on circulant graphs $G=C_m(S)$ where $m=2^n$ for $n\geq 2$. Moreover, it is well know that $G=C_m(S)$ with $S=\{ s_1,s_2,...,s_k\}$ is a connected graph if and only if $gcd(s_1,s_2,...,s_k,m)=1$ (see \cite{Boesch}), 
Therefore, $G=C_m(S)$ is connected if and only if there is an odd integer $l=s_i\in S$ for $i\in [1,k]$.

\vspace{.3cm}
Let $m=2^n$ and $S=\{s_1,s_2,\ldots,s_k\}$ be such that the circulant graph $G=C_m(S)$ is connected and there exists an odd integer $l\in [1, \lfloor \frac{m}{2} \rfloor]\setminus S$. Then $\overline{G}=C_m([1,\lfloor \frac{m}{2} \rfloor]\setminus S)$ is connected. By Observation~\ref{complement},  $G=C_m(S)$ is open XOR-magic if and only if $\overline{G}=C_m([1,\lfloor \frac{m}{2} \rfloor]\setminus S)$ is closed XOR-magic. 
\vspace{.3cm}
The observation below illustrates another connection between open XOR-magic and closed XOR-magic circulant graphs. 

 \begin{observation}
      Let $m=2^n$, $n>1$, and  $S=\{s_1,s_2,\ldots,s_k,m/2\}$ be a subset of $\{1,2,\ldots {m/2}\}$. 
     The graph $G_1=C_m(S)$ is open XOR-magic if and only if $G_2=C_m(S')$ is closed XOR-magic for $S'=\{m/2-s_1,m/2-s_2,\ldots,m/2-s_k\}$.\label{openclosed}
 \end{observation}
 \textit{Proof.} 
 Note that $G_1$ is connected since it is open XOR-magic, thus   there exists $s_j\in S$ that is an odd integer, what implies that $m/2-s_j\in S'$ is odd and $G_2$ is connected.  Assume that $\ell$ is an open XOR-magic labeling of $G_1$. Define a labeling $\ell'$ in $G_2$ as $\ell'(x_i)=\ell(x_i)$, then  \begin{equation} \begin{split}w(x_i)=\sum_{u\in N_{G_1}(x_i)}\ell(u)=\\
 \ell(x_{i-s_1})+\ell(x_{i+s_1})+\ldots+\ell(x_{i-s_k})+\ell(x_{i-s_k})+\ell(x_{i+m/2})\\
 =\sum_{u\in N_{G_2}[x_{i+m/2}]}\ell'(u).
\end{split}
\end{equation}

This finishes the proof.~\qed

\vspace{.3cm}
We introduce the following notation: for any graph $G$, the adjacency matrix of the graph $G$ is denoted by $A(G)$. Let us recall that the adjacency matrix of the circulant graph is a circulant matrix. A square matrix $M$ of size $m$ is called \emph{circulant}, if it satisfies
\begin{align*}
    M = \begin{bmatrix}
a_0     & a_1    & a_2      & \ldots & a_{m-2}  & a_{m-1}\\
a_{m-1}     & a_0    & a_1      & \ldots & a_{m-3}  & a_{m-2}\\
a_{m-2}     & a_{m-1}    & a_0      & \ldots & a_{m-4}  & a_{m - 3}\\
\vdots  & \vdots  & \vdots  & \vdots & \ddots & \vdots \\
a_{1}     & a_{2}    & a_3      & \ldots & a_{m-1}  & a_{0}
\end{bmatrix}.
\end{align*} Such a circulant matrix $M$ is denoted by $\operatorname{circ}_m(a_0,a_1,\ldots a_{m-1})$.

The next theorem give us the connection between the determinant of the adjacency matrix $A(G)$ and the existence of an open XOR-magic labeling for the graph $G$.

\begin{theorem}\label{determinantprop}
   Let $G$ be a graph such that there exists an open XOR-magic labeling of the graph $G$. Then, $\det(A(G))=0 \mod 2$. 
\end{theorem}
    
\noindent \textit{Proof.}
    If $G$ admits an open XOR-magic labeling $\ell$, then $|G| = 2^n = m$ for some natural number $n \in \mathbb{N}$. Let $V(G) = \{x_0, x_1 \ldots, x_{m-1}\}$ and $\mu = (0,\ldots,0)$. Moreover, let $\widetilde{\ell}$ be a column vector of the first coordinates of labels of the vertices of $G$, so $\widetilde{\ell} =  (\ell(x_0)_1, \ell(x_1)_1, \ldots, \ell(x_{m-1})_1)^T$. Then, we have
    \begin{align*}
        A(G) \cdot_2 \widetilde{\ell} = \mu_1 \cdot_2 \mathbf{1},
    \end{align*} where $\mathbf{1}$ is a column vector of length $2^n$ such that each element of $\mathbf{1}$ is equal to one. However, the magic constant $\mu$ is equal to $(0, \ldots, 0)$, so $\mu_1 = 0$ and hence
    \begin{align*}
        A(G) \cdot_2 \widetilde{\ell} = 0 \cdot_2 \mathbf{1} = 0 \cdot_2 \widetilde{\ell}.
    \end{align*} Consequently, over $\mathbb{Z}_2$, the vector $\widetilde{l}$ is an eigenvector of the matrix $A(G)$ with the corresponding eigenvalue $0$ and hence, {$\det(A(G) - 0 \cdot I_{m}) = \det(A(G))=  0 \mod 2 $}.~\qed

\vspace{.3cm}
The inverse of  Theorem \ref{determinantprop} is not true. For example, consider the circulant graph $G = C_{2^n}(\{1\})$ (a cycle of order $2^n$) satisfies $|\det(A(G))| = 2$, but $G$ is not an open XOR-magic graph (because a vertex labeled with $(0,\ldots,0)$ would have to have two identical labels in its neighborhood).

\vspace{.3cm}
In order to see the first application of Theorem \ref{determinantprop}, consider the Cayley graph, known as the hypercube, $\mathcal{Q}_n=\mathrm{Cay}((\mathbb{Z}_2)^n,$ $\{e_1,e_2,\ldots,e_n\})$, where $e_i$ is the vector with a $1$ in the $i$th coordinate and $0$ elsewhere. In \cite{ACFSQ} it was proved that for odd $n$ there is no $\Gamma$-distance magic graph  $\mathcal{Q}_n$ (see Theorem 2.1, \cite{ACFSQ}). In other way, Cvetkovic \cite{Cvetkovic} proved that the hypercube $\mathcal{Q}_n$ has eigenvalues $n - 2i$ for $i\in\{0,1,\ldots,n\}$, so by Theorem ~\ref{determinantprop} we obtain immediately that there is no open XOR-magic $\mathcal{Q}_n$ for $n$ odd.



\vspace{.3cm}
In the rest of this section, we will analyze the nonexistence of an open XOR-magic labeling for different types of circulant graphs.


\begin{lemma}
    Let $m = 2^n$, $n > 2$, $r \in [1,\frac{m}{2} - 1]$ and let $S = \{\frac{m}{2} - r, \frac{m}{2} - (r - 1), \ldots, \frac{m}{2} - 2, \frac{m}{2} - 1, \frac{m}{2}\}$. If $G = C_m(S)$, then $|\det(A(G))| = 2r + 1$.
\end{lemma}

\noindent \textit{Proof}. Observe that  since  $\frac{m}{2}\in S$ and $m - \frac{m}{2} = \frac{m}{2}$, the graph $G$ is odd regular  with valency $2 |S| - 1 = 2 (r+1) - 1 = 2r + 1$ and  the adjacency matrix of the graph $G$ is a circulant matrix
\begin{align*}
    A(G) = \operatorname{circ}_{2^n} v,
\end{align*} where $v = (v_0, v_1, \ldots, v_{m-1})$ is a vector of length $2^n$ such that
\begin{align*}
    v_i = \begin{cases}
        1 & \text{if } i \in \{s, s+r, :s \in S\}, \\
        0 & \text{otherwise}.
    \end{cases}
\end{align*}
Then, we can permute the rows of the matrix $A(G)$ to obtain the circulant matrix $M = \operatorname{circ}_{2^n}(1,1,\ldots,1,0,0,\ldots,0)$ with the first $2r + 1$ elements of the first row equal to $1$ and with other elements in the first row equal to $0$. According to \cite{bincircmatdet} (Theorem 2.1), $\det(M)=2r + 1$ and hence, the determinant satisfies $|\det(A(G))| = 2r + 1$.~\qed

\vspace{.3cm}
{The next corollary is an immediate consequence of  Theorem \ref{determinantprop}.}

\begin{corollary}\label{symetrynot}
    Let $m = 2^n$, $n > 2$, $r \in [1,\frac{m}{2} - 1]$ and let $S = \{\frac{m}{2} - r, \frac{m}{2} - (r - 1), \ldots, \frac{m}{2} - 2, \frac{m}{2} - 1, \frac{m}{2}\}$. If $G = C_m(S)$, then $G$ does not admit an open XOR-magic labeling.
\end{corollary}


\vspace{.3cm}
Remember that if $T=qS$ for some integer $q$ such that $gcd(q,m)=1$, the circulant graph $C_m(S)$ is isomorphic to the circulant graph $C_m(T)$, but the inverse is not true (Adam's conjecture), as we show in the following Theorem.

\begin{theorem}\label{2mminus2}
Let $m = 2^n$, $n > 1$, $S = \{1,\frac{m}{2} -  2,\frac{m}{2}\}={\{1, 2^{n-1} - 2,  2^{n-1}\}}$. If $G = C_m(S)$, then $G$ does not admit an open XOR-magic labeling.
\end{theorem}

\noindent \textit{Proof}. We prove that graph $G$ is isomorphic with a graph $H = C_n(\{\frac{m}{2} - 2, \frac{m}{2} - 1, \frac{m}{2}\})$. Let  $f: V(G) \to V(H)$ be a function defined by a formula
\begin{align*}
    f(x_i) = \begin{cases}
        x_i, & \text{if } 2 | i, \\
        x_{(i + 2^{n - 1}) \mod 2^n} & \text{otherwise}
    \end{cases}.
\end{align*}
Obviously, this is one-to-one mapping. We will prove that it is also edge-preserving function. Let us recall that any two vertices $x_i$ and $x_j$ are adjacent in graph $H$ if and only if $|i - j| = \{\frac{m}{2} - 2,  \frac{m}{2} - 1, \frac{m}{2}, \frac{m}{2} + 1, \frac{m}{2} + 2\}$. Let $x_i x_j \in E(G)$. Without loss of generality, we have the following cases: 

\textbf{Case 1.} $j = (i + 1) \mod 2^n$. If $i$ is even, then $j$ is odd and hence, $f(x_i) = x_i$ and $f(x_{(i+1) \mod 2^n}) = x_{(i+1 + 2^{n-1}) \mod 2^n}$. Moreover, we have 
\begin{align*}
    |(i + 1 + 2^{n - 1})\mod 2^n - i|  = 2^{n-1} + 1 = \frac{m}{2} + 1
\end{align*} and hence, $f(x_i)f(x_j) \in E(H)$. If $i$ is odd, then $j$ is even and hence, $f(x_i) = x_{(i + 2^{n-1}) \mod 2^n}$ and $f(x_j) = x_j$. Again, 
\begin{align*}
    |(i + 2^{n-1}) \mod 2^{n} - (i+1) \mod 2^n| = 2^{n-1} - 1 = \frac{m}{2} - 1,
\end{align*} so $f(x_i) f(x_j) \in E(H)$.

\textbf{Case 2.} $j = (i + \frac{m}{2} - 2) \mod 2^n$. However, $|(i - j) \mod 2^{n}| = \frac{m}{2}$ and hence $f(x_i)f(x_j) \in E(H)$.

Proofs for other cases ($j = (i + 2^{n-1}) \mod 2^n$, $j = \left(i + 2^{n} - \left(\frac{m}{2} - 2\right)\right) \allowbreak \mod 2^n$, $j = (i + 2^{n}-1) \mod 2^n$) are analogous. To conclude, the graph $G$ is isomorphic with the graph $H$.

According to Corollary \ref{symetrynot} for $r = 2$,
the graph $H$ does not admit an open XOR-magic labeling and hence, $G$ is also not $(\mathbb{Z}_2)^n$-distance magic.~\qed

\vspace{.3cm}
For the purposes of further classes of circulant graphs, we will introduce the concept of the Smith normal form of a matrix. For a square integer matrix $M$ of size $n$, a \emph{Smith decomposition} of a matrix $M$ is defined by three matrices $L, R, S$ such that:
\begin{itemize}
    \item $L, R$ are integer, unimodular (the determinant is equal to $1$ or $-1$) invertible matrices of size $n$,
    \item $S$ is a diagonal matrix of size $n$ with the diagonal $(d_0, d_1, \ldots, d_{n-1}) \in \mathbb{N}^n$ such that $d_i | d_{i+1}$ for each $i = 0,1,\ldots,n-2$,
    \item $LMR = S$.
\end{itemize} 
The matrix $S$ is called \emph{Smith normal form} (shortly SNF) of the matrix $M$ and is denoted by $\operatorname{SNF}(M)$. The matrix $\operatorname{SNF}(M)$ is unique. Moreover, $\det(\operatorname{SNF}(M)) \allowbreak = |\det(M)|$. Consequently, for any graph $G$, $\det(\operatorname{SNF}(A(G))) \allowbreak = |\det(A(G))|$ and hence, this form is 
useful for our algebraic approach. Smith normal form of adjacency matrices for various circulant graphs have been studied by Williams (see \cite{williams}). Using the Smith normal forms obtained in this paper and our algebraic approach, we now discuss the problem of the existence of open XOR-magic labelings for some classes of circulant graphs. Diagonal matrix of size $n$ with the diagonal $(x_0,x_1,\ldots,x_{n-1})$ will be denoted by $\operatorname{diag}_n(x_0, x_1, \ldots, x_{n-1})$.

\noindent \textbf{Complement of cycle and of power of cycle}. By $C_m^{(r)}$ we denote a special circulant graph, called \textit{$r$-th power of cycle of length $m$} (where $r \in \mathbb{N}_+$, $r \leq \lfloor \frac{m}{2} \rfloor)$, defined as $C_m^{(r)} = C_m(\{1,2,\ldots,r\})$. Notice that the complement of this graph satisfies $\overline{C_m^{(r)}} = C_m(\{r+1, r+2, \ldots, \lfloor \frac{m}{2} \rfloor\})$ and that $C^{(1)}_m$ is just a cycle.

The adjacency matrix of the complement of a cycle satisfies:

\begin{align*}
    \operatorname{SNF}(A(\overline{C_{n}})) = \begin{cases}
        \operatorname{diag}_n(1,
        \ldots,1,\frac{n-3}{3},0,0) & \text{if } n \mod 3 = 0, \\
        \operatorname{diag}_n(1,
        \ldots,1,n-3) & \text{otherwise} \\
    \end{cases} .
\end{align*} Consequently, if $n = 2^m$, then $|\det(A(\overline{C_n})| = n-3=2^m-3$, which is odd. It gives us the following corollary.

\begin{corollary}
    Let $n \in \mathbb{N}_+$, $n \geq 3$. Then, $\overline{C_n}$ is not an open XOR-magic graph.
\end{corollary}

For the complement of the power of the cycle, Smith normal form of adjacency matrix of this graph satisfies

\begin{align*}
    \operatorname{SNF}(A(\overline{C_{n}^{(r)}})) =
        \operatorname{diag}_n (1,\ldots,1,\frac{n-2r-1}{GCD(n,2r+1)},\underbrace{0, \ldots, 0}_{GCD(n, 2r + 1) - 1}).
\end{align*}

If $n = 2^m$, then $\operatorname{GCD}(n,2r+1) = 1$ and hence, the determinant of matrix $A(\overline{C_n^{(r)}})$ satisfies $|\det(A(\overline{C_n^{(r)}})| = n - 2r - 1 = 2^m - 2r - 1$, which is odd.

\begin{corollary}
    Let $n \in \mathbb{N}_+$, $n \geq 3$ and let $r \in [1,\lfloor \frac{n}{2} \rfloor - 1].$ Then, $\overline{C_n^{(r)}}$ is not an open XOR-magic graph.
\end{corollary}

\noindent \textbf{Andrásfai and Doob graphs}. Let $r \in \mathbb{N}_+$, $r \geq 2$. An \emph{Andrásfai graph} \cite{and} is the circulant graph defined as $\operatorname{And}(r) = C_{3r-1}(\{3k+1: 0 \leq k \leq \lfloor \frac{3r - 3}{2} \rfloor\})$. \emph{Doob graph} is a kind of generalization of a graph $\operatorname{And}(r)$. Let $r, t \in \mathbb{N}_+$ such that $r \geq 2$. A \emph{Doob graph} $G(r,t)$ \cite{doob} is a circulant graph of a form $G(r,t) = C_n(\{kt+1 : 0 \leq k \leq \lfloor \frac{(r-1)t}{2} \rfloor \})$, where $n = (r-1)t + 2$. Notice that $\operatorname{And}(r) = G(r,3)$. Smith normal form of adjacency matrix of Doob graph $G(r,t)$ is given here:
\begin{multline*}
    \operatorname{SNF}(A(G(r,t))) =  \\ \begin{cases}
\operatorname{diag}_{(r-1)t+2} \left( 
1, \ldots, 1, 
\frac{r}{GCD(r,t-2)}, 
\underbrace{0, \ldots, 0}_{GCD(r,t-2)-1}
\right) & \text{if } t \text{ odd}, \\
\operatorname{diag}_{(r-1)t+2} \left( 
1, \ldots, 1, 
\frac{r}{GCD(r,\frac{t-2}{2})-1}, \frac{r}{GCD(r,\frac{t-2}{2})-1}, 
\underbrace{0, \ldots, 0}_{GCD(2r, t-2)-2}
\right) & \text{if } t \text{ even}. \\
    \end{cases}
\end{multline*} If the relevant parameters are relatively prime, then the determinant of the matrix $A(G(r,t))$ satisfies $|\det(A(G(r,t)))| = r$ or $|\det(A(G(r,t))| \allowbreak = r^2$. If $r$ is odd, we obtain the following results.

\begin{corollary}
    Let $r ,t \in \mathbb{N}_+$ be odd natural numbers such that $r \geq 2$ and $\operatorname{GCD}(r,t-2)=1$. Then, the graph $G(r,t)$ is not an open XOR-magic graph.
\end{corollary}

\begin{corollary}
    Let $r \in \mathbb{N}_+$ be odd natural number such that $r \geq 2$ and let $t \in \mathbb{N}_+$ be even natural number. Moreover, let $\operatorname{GCD}(2r,t-2) = 2$. Then, the graph $G(r,t)$ is not an  open XOR-magic graph.
\end{corollary}

\begin{corollary}
    Let $k \in \mathbb{N}_+$. Then, the Andrásfasi graph $\operatorname{And}(k)$ is not an open XOR-magic graph.
\end{corollary}

\noindent \textbf{Möbius ladder graph}. Let $n \in \mathbb{N}_+$ be an even natural number. The \emph{Möbius ladder graph} (defined in \cite{mobius}) is the circulant graph $L_n=C_n(\{1,\frac{n}{2}\})$. If we analyze the existence of open XOR-magic labelings, then we can focus on the ladder graph of order $2^m$. Clearly, $2^m \mod 6 \in \{\pm 1, \pm 2\}$. Smith normal form for ladder graph $L_n$ with such parameters is as follows:
\begin{align*}
    SNF(A(L_n)) = \begin{cases}
        \operatorname{diag}_n(1,1,\ldots,3,3) & \text{if } \frac{n}{2} \mod 6 \in \{1, -1\}, \\
        \operatorname{diag}_n(1,1,\ldots,1,3) & \text{if } \frac{n}{2} \mod 6 \in \{2, -2\}, \\
    \end{cases}
\end{align*}
and hence, $|\det(A(L_{2^m}))| \in \{3,9\}$.
\begin{corollary}
    Let $n \in \mathbb{N}_+$. Then, Möbius ladder graph $L_n$ is not an open XOR-magic graph.
\end{corollary}

\section{Odd regular open XOR-magic graphs and even regular closed XOR-magic graphs}\label{odd_even}
 In this section, we prove Theorem \ref{th:C1P}, thereby giving a positive answer to Open Problem \ref{XORopenproblem} (mentioned in the Introduction) posed by Batal \cite{Batal}. The proof proceeds by induction on  $n\geq 4$, the number of vertices, and relies on known results on graph products, together with Theorems \ref{openxomagicevenregular} and \ref{closedxomagicevenregular}, which form the basis of our inductive argument.

To establish the existence of the graphs related to the next two theorems, we develop several Mixed Integer Linear Programming (MILP) models, which we present in the final part of this section.

\begin{restatable}{theorem}{CoPSThm}
  \label{openxomagicevenregular}
For any $n$ and $k$ such that
\begin{enumerate}
    \item $n=4$ and $k\in \{5,7,9,11\}$, 
    \item $n=5$ and $k\in \{5,7,9,11,13,15,27\}$,
    \item $n = 6$ and $k = 9$,
    \item $n = 7$ and $k = 13$,
\end{enumerate}
there exists an open XOR-magic $k$-regular graph of power $n$.
\end{restatable}

\begin{restatable}{theorem}{CoPCThm}
\label{closedxomagicevenregular}
For any $n$ and $k$ such that
\begin{enumerate}
    \item $n=4$ and $k\in \{4,6,8,10\}$,
    \item $n = 5$ and $k \in \{26,24, 22, 20, 18, 16, 4\}$,
    \item $n = 6$ and $k = 54$,
    \item $n = 7$ and $k = 114$,
\end{enumerate}
there exists a closed XOR-magic $k$-regular graph of power $n$.
\end{restatable}

Let us recall two out of four standard graph products (see \cite{IK}). Both the \emph{Cartesian product} $G\square H$ and the \emph{strong product} $G\boxtimes H$ are graphs with vertex set $V(G)\times V(H)$. Two vertices $(g,h)$ and $(g^{\prime },h^{\prime })$ are adjacent in:

\begin{itemize}
\item $G\square H$ if and only if either $g=g^{\prime }$ and $h$ is adjacent
with $h^{\prime }$ in $H$, or $h=h^{\prime }$ and $g$ is adjacent with $
g^{\prime }$ in $G$;


\item $G\boxtimes H$ if and only if  $g=g^{\prime }$ and $h$ is adjacent
with $h^{\prime }$ in $H$, or $h=h^{\prime }$ and $g$ is adjacent with $
g^{\prime }$ in $G$, or $g$ is adjacent with $g^{\prime }$ in $G$
and $h$ is adjacent with $h^{\prime }$ in $H$.

\end{itemize}
From the above definition, it follows that if
 $G$ is a $k$-regular graph  and $H$ is an $l$-regular graph, then $G\square H$ is $(k+l)$-regular and $G\boxtimes H$ is $(kl+k+l)$-regular. 

  \begin{lemma}[\cite{IK}]  Let $G$ and $H$ be simple graphs, then
 \begin{itemize}
    
     \item $G\square H$ is connected if and only if each of its factors is connected;

\item $G\boxtimes H$ is connected if and only if each of its factors is connected.
 \end{itemize} 
 \end{lemma}

The following was proved:
\begin{corollary}[\cite{Batal}]\label{cartesian} Let $G$ and $H$ be simple graphs, then
\begin{itemize}
    \item $G\square H$ is  an even-regular closed XOR-magic graph if both factors are  even-regular closed XOR-magic graphs or  odd-regular open XOR-magic graphs;
      \item $G\square H$ is an odd-regular open XOR-magic graph if one factor is an
odd-regular open XOR-magic graph, and the other is an even-regular closed XOR-magic graph;
\item $G\boxtimes H$ is  an even-regular closed XOR-magic graph if both factors are  even-regular closed XOR-magic graphs;
\end{itemize}
    
\end{corollary}

From the above corollary, we obtain the following main result.

\CoPthm*

\textit{Proof.}
We will proceed by induction. For $n\in\{4,5,6,7\}$, there exists an even closed XOR-magic graph and an odd open XOR-magic graph of power $n$ by Theorem~\ref{closedxomagicevenregular} and \ref{openxomagicevenregular} (explained in detail later in this section), respectively. 
 Let $H_1$ be an odd open XOR-magic graph of power $4$ and  $H_2$ be an even closed XOR-magic graph of power $4$. Let $n\geq 8$,  there exists an even closed XOR-magic graph $G_1$ of power $n-4$ and an odd open XOR-magic graph $G_2$ by the induction hypothesis. 
 Thus, the graph $G'=G_1\square H_2$ is an even closed XOR-magic graph, whereas $G'=G_1\square H_1$ is an odd open XOR-magic graph of power $n$ by Corollary~\ref{cartesian}.~\qed

\textbf{Finding the base graphs}. In this part of the section, we describe how the graphs from Theorem \ref{openxomagicevenregular} and Theorem \ref{closedxomagicevenregular} were obtained. To obtain them, we have developed some Mixed Integer Linear Programming (MILP) models.

Firstly, let us focus on the existence of odd regular open XOR-magic graphs. 

\vspace{.3cm}

The first model is a formulation of the problem of the existence of a $d$-regular graph with $2^n$ vertices such that there exists an open XOR-magic labeling. In details, let $d \in \mathbb{N}_+$ be fixed odd number, let $n \in \mathbb{N}_+$ be fixed power in the number of vertices and let $\ell: [1, 2^n] \to (\mathbb{Z}_2)^n$ be a fixed bijection (we associate a set of vertices with the set $[1, 2^n]$, for simplifying the notation, let $V = [1, 2^n]$). Then, we introduce the variables in our model:
\begin{itemize}
    \item $e_{u,v} \in \{0,1\}$ for all pairs $(u, v) \in V^2$ such that $u \neq v$ -- binary decision variable: $1$ if $uv \in E$ and $0$ otherwise.
    \item $k_{u,i} \in \mathbb{Z}$ for each node $u \in V$ and each index $i \in [n]$ -- this variable allow us to express the sum of labels of the neighbors of the node $u$ at index $i$. 
\end{itemize}

Then, we have to introduce the constraints. Firstly, our graph is simple, undirected graph, so the adjacency matrix of searched graph must be symmetric. Consequently, we introduce constraints: 
\begin{align}\label{edgesConstraints1}
    e_{u,v} = e_{v,u} \text{ for each pair } (u, v) \in V^2, u \neq v.
\end{align} Our graph must be $d$-regular, so we have another constraints:
\begin{align}\label{edgesConstraints2}
    \sum_{v \in [2^m] - \{u\}} e_{u,v} = d \quad \text{for each node } u \in V.
\end{align}

Then, for each node $v \in [2^n]$, sum of labels over open neighborhood of $v$ at index $i$ must be equal to $0$. Equivalently, if we add labels in the sense of adding natural numbers at each index, this sum must be equal to $(2k_{v,1}, 2k_{v,2}, \ldots, 2k_{v,n})$. It gives us another conditions:
\begin{align}\label{pairity}
    2 k_{v,i} = \sum_{u \in [2^n] - \{v\}} (\ell(u)_i \cdot e_{u,v}) \quad \text{for each node } v \in [2^n] \text{ and each index } i \in [n].
\end{align}
Finally, all we need to do is set an objective function in our model. 
It just needs to be a constant function. Formally, this model is presented in Figure \ref{modelclassic}.

\begin{figure}
    \begin{equation*}
\begin{aligned}
\text{min} \quad & \text{\texttt{const}} \\
\text{s.t.} \quad 
& e_{u,v} = e_{v,u} & \text{for each } (u,v) \in V^2, u \neq v, \\
& \sum_{v \in V - \{u\}} e_{u,v} = d & \text{for each } u \in V, \\
& 2k_{v,i} = \sum_{u \in V - \{v\}} (\ell(u)_i \cdot e_{u,v}) & \text{for each } v \in V \text{ and each } i=1,2,\ldots, n, \\
& e_{u,v} \in \{0,1\} & \text{for each } (u, v) \in V^2, u \neq v, \\
& k_{v, i} \in \mathbb{Z} & \text{for each } v \in V \text{ and each } i =1,2,\ldots,n
\end{aligned}
\end{equation*}
    \caption{First MILP model for finding odd regular open XOR-magic graphs (parameters: $n \in \mathbb{N}_+, V = [1,2^n], d \in \mathbb{N}_+$ - odd number, $\ell: V \to (\mathbb{Z}_2)^n$ - fixed bijection)}
    \label{modelclassic}
\end{figure}

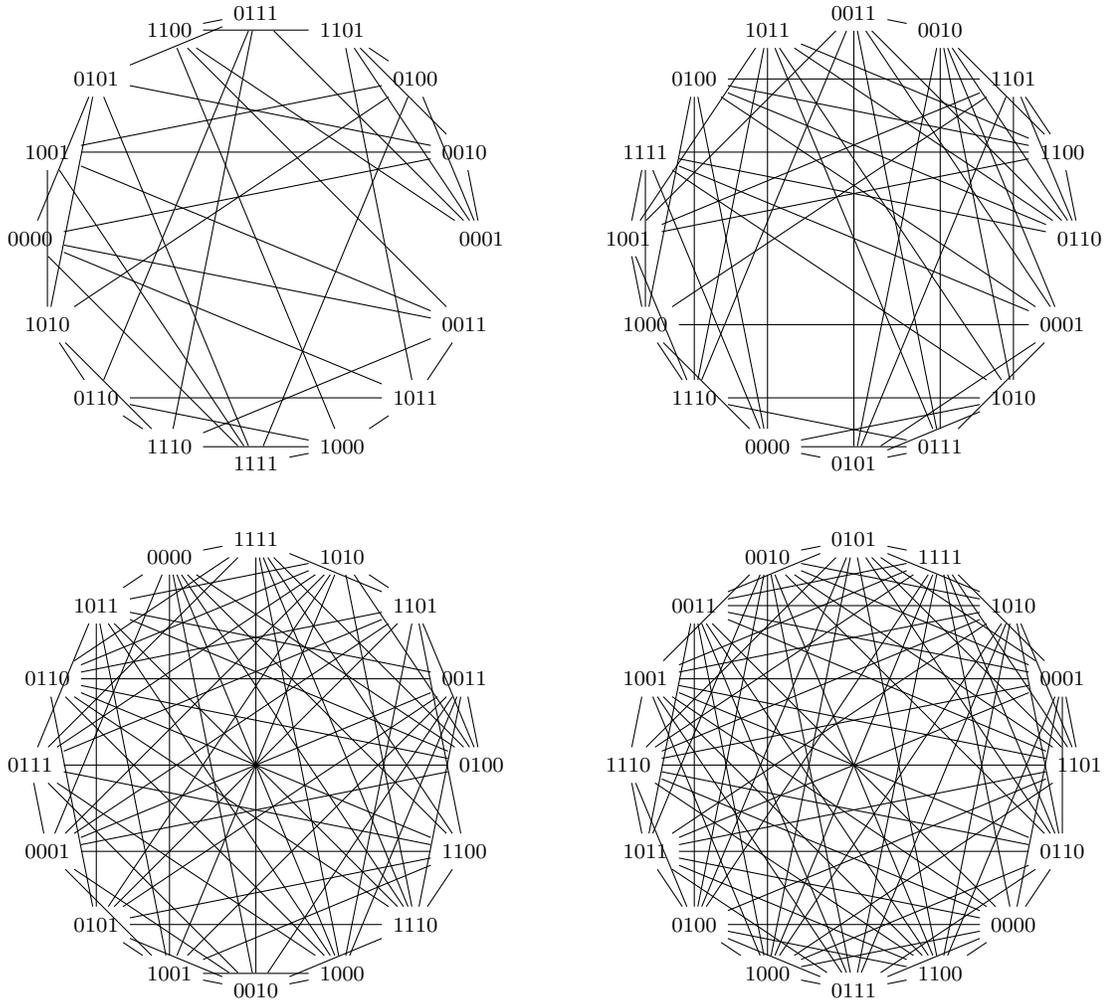
\begin{figure}[!h]
    \centering
    \begin{minipage}{0.42\textwidth}
    \begin{tikzpicture}
      \draw
        (0.0:3) node (0){\scriptsize\scriptsize $0001$}
        (22.5:3) node (1){\scriptsize$0010$}
        (45.0:3) node (3){\scriptsize$0100$}
        (67.5:3) node (7){\scriptsize$1101$}
        (90.0:3.) node (8){\scriptsize$0111$}
        (112.5:3) node (11){\scriptsize$1100$}
        (135.0:3) node (4){\scriptsize$0101$}
        (157.5:3) node (9){\scriptsize$1001$}
        (180.0:3.) node (12){\scriptsize$0000$}
        (202.5:3) node (2){\scriptsize$1010$}
        (225.0:3) node (5){\scriptsize$0110$}
        (247.5:3) node (15){\scriptsize$1110$}
        (270.0:3) node (14){\scriptsize$1111$}
        (292.5:3) node (6){\scriptsize$1000$}
        (315.0:3) node (10){\scriptsize$1011$}
        (337.5:3) node (13){\scriptsize$0011$};
      \begin{scope}[-]
        \draw (0) to (1);
        \draw (0) to (3);
        \draw (0) to (7);
        \draw (0) to (8);
        \draw (0) to (11);
        \draw (1) to (4);
        \draw (1) to (7);
        \draw (1) to (9);
        \draw (1) to (12);
        \draw (3) to (2);
        \draw (3) to (7);
        \draw (3) to (9);
        \draw (3) to (14);
        \draw (7) to (10);
        \draw (7) to (11);
        \draw (8) to (4);
        \draw (8) to (5);
        \draw (8) to (11);
        \draw (8) to (15);
        \draw (11) to (6);
        \draw (11) to (13);
        \draw (4) to (2);
        \draw (4) to (12);
        \draw (4) to (14);
        \draw (9) to (2);
        \draw (9) to (13);
        \draw (9) to (14);
        \draw (12) to (10);
        \draw (12) to (13);
        \draw (12) to (14);
        \draw (2) to (5);
        \draw (2) to (15);
        \draw (5) to (6);
        \draw (5) to (10);
        \draw (5) to (15);
        \draw (15) to (6);
        \draw (15) to (13);
        \draw (14) to (6);
        \draw (6) to (10);
        \draw (10) to (13);
      \end{scope}
    \end{tikzpicture}
    \end{minipage}
    \hfill
  \begin{minipage}{0.42\textwidth}
    \centering
      \begin{tikzpicture}
      \draw
        (0.0:3) node (0){\scriptsize $0110$}
        (22.5:3) node (1){\scriptsize $1100$}
        (45.0:3) node (4){\scriptsize $1101$}
        (67.5:3) node (7){\scriptsize $0010$}
        (90.0:3) node (8){\scriptsize $0011$}
        (112.5:3) node (10){\scriptsize $1011$}
        (135.0:3) node (11){\scriptsize $0100$}
        (157.5:3) node (15){\scriptsize $1111$}
        (180.0:3) node (2){\scriptsize $1001$}
        (202.5:3) node (5){\scriptsize $1000$}
        (225.0:3) node (12){\scriptsize $1110$}
        (247.5:3) node (3){\scriptsize $0000$}
        (270.0:3) node (6){\scriptsize $0101$}
        (292.5:3) node (9){\scriptsize $0111$}
        (315.0:3) node (13){\scriptsize $1010$}
        (337.5:3) node (14){\scriptsize $0001$};
      \begin{scope}[-]
        \draw (0) to (1);
        \draw (0) to (4);
        \draw (0) to (7);
        \draw (0) to (8);
        \draw (0) to (10);
        \draw (0) to (11);
        \draw (0) to (15);
        \draw (1) to (2);
        \draw (1) to (4);
        \draw (1) to (7);
        \draw (1) to (10);
        \draw (1) to (11);
        \draw (1) to (15);
        \draw (4) to (2);
        \draw (4) to (5);
        \draw (4) to (6);
        \draw (4) to (11);
        \draw (4) to (13);
        \draw (7) to (6);
        \draw (7) to (8);
        \draw (7) to (9);
        \draw (7) to (13);
        \draw (7) to (14);
        \draw (8) to (2);
        \draw (8) to (6);
        \draw (8) to (9);
        \draw (8) to (12);
        \draw (8) to (14);
        \draw (10) to (2);
        \draw (10) to (3);
        \draw (10) to (9);
        \draw (10) to (12);
        \draw (10) to (13);
        \draw (11) to (3);
        \draw (11) to (5);
        \draw (11) to (12);
        \draw (11) to (14);
        \draw (15) to (2);
        \draw (15) to (3);
        \draw (15) to (5);
        \draw (15) to (13);
        \draw (15) to (14);
        \draw (2) to (5);
        \draw (2) to (12);
        \draw (5) to (3);
        \draw (5) to (12);
        \draw (5) to (14);
        \draw (12) to (9);
        \draw (12) to (13);
        \draw (3) to (6);
        \draw (3) to (9);
        \draw (3) to (13);
        \draw (6) to (9);
        \draw (6) to (13);
        \draw (6) to (14);
        \draw (9) to (14);
      \end{scope}
    \end{tikzpicture}
  \end{minipage}
  \vspace{0.5cm} 
  \begin{minipage}{0.42\textwidth}
      \begin{tikzpicture}
      \draw
        (0.0:3) node (0){\scriptsize$0100$}
        (22.5:3) node (1){\scriptsize$0011$}
        (45.0:3) node (2){\scriptsize$1101$}
        (67.5:3) node (3){\scriptsize$1010$}
        (90.0:3) node (6){\scriptsize$1111$}
        (112.5:3) node (7){\scriptsize$0000$}
        (135.0:3) node (8){\scriptsize$1011$}
        (157.5:3) node (11){\scriptsize$0110$}
        (180.0:3) node (12){\scriptsize$0111$}
        (202.5:3) node (13){\scriptsize$0001$}
        (225.0:3) node (4){\scriptsize$0101$}
        (247.5:3) node (9){\scriptsize$1001$}
        (270.0:3) node (10){\scriptsize$0010$}
        (292.5:3) node (14){\scriptsize$1000$}
        (315.0:3) node (15){\scriptsize$1110$}
        (337.5:3) node (5){\scriptsize$1100$};
      \begin{scope}[-]
        \draw (0) to (1);
        \draw (0) to (2);
        \draw (0) to (3);
        \draw (0) to (6);
        \draw (0) to (7);
        \draw (0) to (8);
        \draw (0) to (11);
        \draw (0) to (12);
        \draw (0) to (13);
        \draw (1) to (4);
        \draw (1) to (8);
        \draw (1) to (9);
        \draw (1) to (10);
        \draw (1) to (11);
        \draw (1) to (13);
        \draw (1) to (14);
        \draw (1) to (15);
        \draw (2) to (3);
        \draw (2) to (4);
        \draw (2) to (5);
        \draw (2) to (6);
        \draw (2) to (11);
        \draw (2) to (12);
        \draw (2) to (13);
        \draw (2) to (14);
        \draw (3) to (4);
        \draw (3) to (8);
        \draw (3) to (9);
        \draw (3) to (11);
        \draw (3) to (12);
        \draw (3) to (13);
        \draw (3) to (15);
        \draw (6) to (5);
        \draw (6) to (7);
        \draw (6) to (10);
        \draw (6) to (11);
        \draw (6) to (12);
        \draw (6) to (14);
        \draw (6) to (15);
        \draw (7) to (5);
        \draw (7) to (8);
        \draw (7) to (9);
        \draw (7) to (10);
        \draw (7) to (13);
        \draw (7) to (14);
        \draw (7) to (15);
        \draw (8) to (4);
        \draw (8) to (9);
        \draw (8) to (12);
        \draw (8) to (14);
        \draw (8) to (15);
        \draw (11) to (4);
        \draw (11) to (5);
        \draw (11) to (14);
        \draw (11) to (15);
        \draw (12) to (5);
        \draw (12) to (10);
        \draw (12) to (13);
        \draw (12) to (14);
        \draw (13) to (5);
        \draw (13) to (9);
        \draw (13) to (10);
        \draw (4) to (5);
        \draw (4) to (9);
        \draw (4) to (10);
        \draw (4) to (15);
        \draw (9) to (5);
        \draw (9) to (10);
        \draw (9) to (14);
        \draw (10) to (14);
        \draw (10) to (15);
        \draw (15) to (5);
      \end{scope}
    \end{tikzpicture}
    \end{minipage}
    \hfill
  \begin{minipage}{0.42\textwidth}
    \centering
      \begin{tikzpicture}
      \draw
        (0.0:3) node (0){\scriptsize$1101$}
        (22.5:3) node (2){\scriptsize$0001$}
        (45.0:3) node (4){\scriptsize$1010$}
        (67.5:3) node (7){\scriptsize$1111$}
        (90.0:3) node (8){\scriptsize$0101$}
        (112.5:3) node (9){\scriptsize$0010$}
        (135.0:3) node (10){\scriptsize$0011$}
        (157.5:3) node (11){\scriptsize$1001$}
        (180.0:3) node (12){\scriptsize$1110$}
        (202.5:3) node (13){\scriptsize$1011$}
        (225.0:3) node (14){\scriptsize$0100$}
        (247.5:3) node (15){\scriptsize$1000$}
        (270.0:3) node (1){\scriptsize$0111$}
        (292.5:3) node (3){\scriptsize$1100$}
        (315.0:3) node (5){\scriptsize$0000$}
        (337.5:3) node (6){\scriptsize$0110$};
      \begin{scope}[-]
        \draw (0) to (2);
        \draw (0) to (4);
        \draw (0) to (7);
        \draw (0) to (8);
        \draw (0) to (9);
        \draw (0) to (10);
        \draw (0) to (11);
        \draw (0) to (12);
        \draw (0) to (13);
        \draw (0) to (14);
        \draw (0) to (15);
        \draw (2) to (3);
        \draw (2) to (5);
        \draw (2) to (6);
        \draw (2) to (7);
        \draw (2) to (8);
        \draw (2) to (9);
        \draw (2) to (10);
        \draw (2) to (11);
        \draw (2) to (12);
        \draw (2) to (13);
        \draw (4) to (1);
        \draw (4) to (3);
        \draw (4) to (6);
        \draw (4) to (8);
        \draw (4) to (9);
        \draw (4) to (10);
        \draw (4) to (11);
        \draw (4) to (12);
        \draw (4) to (13);
        \draw (4) to (15);
        \draw (7) to (1);
        \draw (7) to (5);
        \draw (7) to (6);
        \draw (7) to (8);
        \draw (7) to (10);
        \draw (7) to (11);
        \draw (7) to (12);
        \draw (7) to (14);
        \draw (7) to (15);
        \draw (8) to (3);
        \draw (8) to (6);
        \draw (8) to (9);
        \draw (8) to (10);
        \draw (8) to (12);
        \draw (8) to (14);
        \draw (8) to (15);
        \draw (9) to (1);
        \draw (9) to (3);
        \draw (9) to (6);
        \draw (9) to (11);
        \draw (9) to (13);
        \draw (9) to (14);
        \draw (9) to (15);
        \draw (10) to (1);
        \draw (10) to (3);
        \draw (10) to (5);
        \draw (10) to (13);
        \draw (10) to (14);
        \draw (10) to (15);
        \draw (11) to (1);
        \draw (11) to (3);
        \draw (11) to (5);
        \draw (11) to (6);
        \draw (11) to (12);
        \draw (11) to (15);
        \draw (12) to (3);
        \draw (12) to (5);
        \draw (12) to (6);
        \draw (12) to (13);
        \draw (12) to (14);
        \draw (13) to (1);
        \draw (13) to (3);
        \draw (13) to (5);
        \draw (13) to (6);
        \draw (13) to (14);
        \draw (14) to (1);
        \draw (14) to (3);
        \draw (14) to (5);
        \draw (14) to (15);
        \draw (15) to (1);
        \draw (15) to (5);
        \draw (15) to (6);
        \draw (1) to (3);
        \draw (1) to (5);
        \draw (1) to (6);
        \draw (3) to (5);
        \draw (5) to (6);
      \end{scope}
    \end{tikzpicture}
  \end{minipage}
    \caption{Open XOR-Magic $d$-regular graphs of order $2^4$ for $d \in \{5,7,9,11\}$.}
    \label{openXORoddregularfig}
\end{figure}

Notice that if this model is satisfiable, then there exists an odd-regular graph of order $2^n$ such that this graph admits an open XOR-magic labeling. However, this graph may not be an open XOR-magic graph, because this graph can be disconnected. Still, if the returned graph is connected, we will find an odd regular open XOR-magic graph. It occurs that this model is satisfiable for some parameters $n$ and $d$ and the generated graphs are connected. In particular, some odd regular open XOR-magic graphs of order $16$ are depicted in Figure \ref{openXORoddregularfig}.

The presented model is simply an expression of the existence of an open XOR-magic labeling in terms of linear and integer programming. We want to introduce another, more efficient model to this problem. Our main idea is to reduce the number of constraints without making them too long. To obtain this, we will represent sequences as integers, and the sequences will correspond to the representations of these integers in the appropriate numeral systems. More general, we introduce a strategy of encoding of binary sequences as integers.

\begin{definition}
    Let $r \in \mathbb{N}_+$ and let $s = (s_1, s_2, \ldots, s_r) \in \mathbb{N}_+^r$ be a sequence of positive numbers such that $s_1 > s_2 >\ldots > s_r$. Moreover, let $x = (x_1, x_2, \ldots, x_r) \in \mathbb{Z}_2^r$ be a binary sequence of length $r$. Then, \emph{$s$-code} of a sequence $x$ is defined as
    \begin{align*}
        \operatorname{enc}_{s}(x) = \sum_{i=1}^{r} s_i \cdot x_i.
    \end{align*}
\end{definition}

For example, if $x = (1,0,1,1)$ and $s = (9^3, 9^2, 9^1, 9^0)$, then $\operatorname{enc}_s(x) = 1 \cdot 9^3 + 0 \cdot 9^2 + 1 \cdot 9^1 + 1 \cdot 9^0 = 739$. Notice that in this case, sequence $x$ is just a base-$9$ representation of a number $739$.

\begin{lemma}\label{encodingLemma}
    Let:
    \begin{itemize}
        \item $M, k, r \in \mathbb{N}_+$, $y \in \mathbb{N}$,
        \item $S_r = (M^{r - 1 + y}, M^{r-2 + y}, \ldots, M^{y})$,
        \item $k_1, k_2, \ldots, k_r \in \mathbb{N}$ such that $0 \leq 2 k_i < M$ for each $i = 1, 2, \ldots r$.
    \end{itemize} Moreover, let $\mathcal{A} = \{a_1, a_2, \ldots, a_k\}$ be a set of binary sequences of length $r$ such that
    \begin{itemize}
        \item $M > \displaystyle \max_{i=1,2,\ldots,r} |\{j : (a_j)_i = 1\}|$,
        \item $\sum_{j=1}^{k} \operatorname{enc}_{S_r}(a_j) = \sum_{i=1}^{r} 2 k_i \cdot M^{y + r - i}$.
    \end{itemize} Then, we have
    \begin{align*}
        \sum_{j=1}^{k} a_j = (0,0,\ldots,0).
    \end{align*}
    and for each $i=1,2, \ldots, r$ we have $|\{j : (a_j)_i = 1| =  2k_i$.
\end{lemma}

\textit{Proof.} Without loss of generality, assume that $y = 0$. Then, we define a sequence $R = (R_1, R_2, \ldots, R_r) \in \mathbb{N}^{k}$ such that $R_i = (a_1)_i +_{\mathbb{N}} (a_2)_i +_{\mathbb{N}} \ldots +_{\mathbb{N}} (a_k)_i$ for each $i = 1,2, \ldots, r$. 

Notice that each element $R_i$ of this sequence satisfies 
\begin{align*}
    0 \leq R_i \leq \displaystyle \max_{w = 1,2,\ldots,r} |\{j : (a_j)_w = 1\}| < M.
\end{align*} As a consequence, we can interpret the sequence $R$ as base-$M$ representation of some natural number $\gamma$. This number $\gamma$ satisfies:

\begin{multline*}
    \gamma = \sum_{i=1}^{r} R_i \cdot M^{k-i} = \sum_{i=1}^{r} M^{k-i} \left( \sum_{j=1}^{r} (a_j)_i \right) = \\ \sum_{j=1}^{r} \left( \sum_{i=1}^{k} M^{k-i} \cdot (a_j)_i \right) = \sum_{j=1}^{k} \operatorname{enc}_{S_r}(a_j).
\end{multline*}

From the assumption we know that $\sum_{j=1}^{k} \operatorname{enc}_{S_r}(a_j) = \sum_{i=1}^{r} 2k_i M^{k-i}$ and $0 \leq 2 k_i < M$ for each $i = 1,2, \ldots, r$. It gives us that a sequence $\mathcal{K} = (2 k_1, 2k_2, \ldots, 2k_r)$ is also base-$M$ representation of the same natural number $\gamma$. Base-$M$ representation of a natural number is unique. Consequently, we have $R = \mathcal{K}$. 

Moreover, for each $i = 1,2, \ldots, r$, the following equality holds
\begin{align*}
    \sum_{j=1}^{k} (a_j)_i = 2k_i
\end{align*} and hence, $a_1 + a_2 + \ldots + a_k = (0,0,\ldots,0)$.~\qed

\begin{remark}\label{remarkOdd}
    Notice that if $M$ in Lemma \ref{encodingLemma} is odd, then assumption $M > \displaystyle \max_{i=1,2,\ldots, r} |\{ j : (a_j)_i = 1|$ can be replaced by assumption $$M \geq \displaystyle \max_{i=1,2,\ldots, r} |\{ j : (a_j)_i = 1|,$$ because equality in this inequality will never occur under the remaining assumptions. Assuming that this equality holds, we will get a sequence representing the natural number $\gamma = \sum_{j=1}^{k} \operatorname{enc}_{S_r}(a_j)$ in base-$M$ system such that at least one element in this representation is odd. However, another representation of the same number $\gamma$ in this system is the sequence $\mathcal{K} = (2k_1, 2k_2, \ldots, 2k_r)$ and the representation is unique, a contradiction.
\end{remark}

Let us now present second MILP using this encoding and the above observation. 

\vspace{.3cm}

Again, let $d \in \mathbb{N}_+$ be fixed odd number, let $n \in \mathbb{N}_+$ be fixed power in the number of vertices and let $\ell: [1, 2^n] \to (\mathbb{Z}_2)^n$ be a fixed bijection (we associate a set of vertices with the set $[1, 2^n]$). We will split each label into smaller sublabels and then, encode each sublabel as integer. Formally, let $t \in \mathbb{N}_+$ be a parameter of the model such that $t \geq 1$ and $t \leq n$ and let $\widetilde{\ell}$ be a function defined for a set $[1, 2^n]$ by a formula
\begin{multline*}
    \widetilde{\ell}(v) =  ( (\ell(v)_1, \ell(v)_2,\ldots, \ell(v)_t), (\ell(v)_{t+1}, \ell(v)_{t+2}, \ldots, \ell(v)_{2t}), \ldots , \\(\ell(v)_{\left(\lceil \frac{n}{t} \rceil - 1\right) t + 1}, \ldots, \ell(v)_n)).
\end{multline*} For example, if $n = 7$, $t = 3$ and $\ell(v) = (1,0,0,1,1,1,0)$ for some $v$, then $\widetilde{\ell}(v) = ((1,0,0),(1,1,1),(0))$. We introduce the same variables $e_{u,v}$ and $k_{u, i}$ as in the previous model and we copy constraints \eqref{edgesConstraints1} and \eqref{edgesConstraints2}. Next, we want to encode each subsequence using our encoding, based on the sequence $(\ldots, d^3, d^2, d)$, while handling the condition from the definition of open XOR-magic labeling. For example, if $d = 7$ and $\ell(v) = (1,0,0,1,1,1,0)$, then $\widetilde{\ell(v)} = ((1,0,0),(1,1,1),(0))$ and we encode this as triple $(1 \cdot 7^3, 1 \cdot 7^2 + 1 \cdot 7^2 + 1 \cdot 7, 0)$. To obtain this, we add constraints $0 \leq 2 k_{v,i} \leq d-1$ for each vertex $v \in V$ and each $i = 1,2, \ldots, n$ and we replace constraints \eqref{pairity} with the following:
\begin{multline*}
    \sum_{i=1}^{t} 2 k_{v, (q-1) \cdot t + i} \cdot d^{t-i+1} = \sum_{u \in [2^n] - \{v\}} (\operatorname{enc}_{S_t}([\widetilde{\ell}(u)]_q) \cdot e_{u,v}) \\ \text{for each node } v \in V \text{ and each } q = 1,2, \ldots, \left\lceil \frac{n}{t} \right\rceil - 1
\end{multline*} 
(where $S_t = (d^t, d^{t-1}, \ldots, d)$) and
\begin{multline*}
    \sum_{i=\left(\lceil \frac{n}{t} \rceil - 1\right)t+1}^{n} 2 k_{v,i} d^{n-i+1} = \sum_{u \in V - \{v\}} (\operatorname{enc}_{S_{end}} ([\widetilde{\ell}(u)]_{\lceil \frac{n}{t} \rceil}) \cdot e_{u, v}) \\ \text{for each node } v \in V
\end{multline*} (where $S_{end} = ( d^{n - \left(\lceil \frac{n}{t} \rceil - 1 \right)t}, \ldots, d^2, d)$).

According to Lemma \ref{encodingLemma} and Remark \ref{remarkOdd}, if these constraints are satisfied, then for each vertex $v \in V$, sum of labels  over open neighborhood of $v$ is equal to $(0,0,\ldots,0)$. Notice that thanks to our encoding, we can 
significantly reduce the number of constraints of the last type: for each vertex $v \in V$, we reduce this number of constraints from $n$  to $\lceil \frac{n}{t} \rceil$. 
By the value of parameter $t$, we can control the length of each of above constraints. Formally, our second MILP model is presented in Figure \ref{secondModelFig}.

\newpage
\begin{figure}
    \begin{equation*}
\begin{aligned}
\text{min} \quad & \text{\texttt{const}} \\
\text{s.t.} \quad 
& e_{u,v} = e_{v,u} & \text{for each } (u,v) \in V^2, u \neq v, \\
& \sum_{v \in V - \{u\}} e_{u,v} = d & \text{for each } u \in V, \\
&  \sum_{i=1}^{t} 2 k_{v, (q-1) \cdot t + i} \cdot d^{t-i+1} =  & \\
& \quad = \sum_{u \in [2^n] - \{v\}} (\operatorname{enc}_{S_t}([\widetilde{\ell}(u)]_q) \cdot e_{u,v}) & \text{for } v \in V, q \in \left\{1, \ldots, \left\lceil \frac{n}{t} \right\rceil -1\right\}, & \\
& \sum_{i=\left(\lceil \frac{n}{t} \rceil - 1 \right) t + 1}^{n} 2 k_{v,i} d^{n-i+1} = & \\ 
& \quad = \sum_{u \in V - \{v\}} (\operatorname{enc}_{S_{end}} ([\widetilde{\ell}(u)]_{\lceil \frac{n}{t} \rceil}) \cdot e_{u, v}) & \text{for each } v \in V \\
& e_{u,v} \in \{0,1\} & \text{for each } (u, v) \in V^2, u \neq v, \\
& k_{v, i} \in \mathbb{Z} & \text{for each } v \in V, i =1,2,\ldots,n, \\
& 0 \leq 2 k_{v, i} \leq d-1 & \text{for each } v \in V, i =1,2,\ldots,n, \\
\end{aligned}
\end{equation*}
    \caption{Second MILP model for finding odd regular open XOR-magic graphs (parameters: $n \in \mathbb{N}_+, t \in \mathbb{N}_+, t \geq 1, t \leq n, V = [1,2^n], d \in \mathbb{N}_+$ - odd number, $\ell: V \to (\mathbb{Z}_2)^n$ - fixed bijection, $S_t = (d^t, d^{t-1}, \ldots, d), S_{end} = ( d^{n - \left(\lceil \frac{n}{t} \rceil - 1 \right)t}, \ldots, d^2, d)$)}
    \label{secondModelFig}
\end{figure}

The second model is definitely more efficient than the first MILP model. With this model, we found an odd regular (in details, $13$-regular) open XOR-magic graph of order $2^7$ with $10$ minutes of computation (the first model did not provide a solution for this order even after days of computation). In details, the graphs we found are listed below.

\CoPSThm*

All odd regular open XOR-magic graphs that we have found can be found at \url{https://github.com/grochowskih/RegularXorMagicGraphs}.

If we go to even regular closed XOR-magic graphs, we can easily adapt the models for odd regular closed XOR-magic graphs. It is enough to define the model parameters in a proper way and modify the constraints to include also the label of a given vertex, not only labels of neighbors. 

We can also utilize odd regular open XOR-magic graphs that we have found. Easily, if $G$ is an odd regular open XOR-magic graph, then the complement of this graph $\overline{G}$ is even regular and $\overline{G}$ admits $(\mathbb{Z}_2)^n$-closed distance magic labeling. Additionally, if $\overline{G}$ is connected, then $\overline{G}$ is even regular closed XOR-magic graph. 

With the above observations, we obtained the following theorem.

\CoPCThm*

\begin{figure}[!h]
    \centering
      \begin{tikzpicture}[scale=1.8]
      \draw
        (0.0:2) node (0){\scriptsize $1101$}
        (22.5:2) node (2){\scriptsize $0001$}
        (45.0:2) node (4){\scriptsize $1010$}
        (67.5:2) node (7){\scriptsize $1111$}
        (90.0:2) node (8){\scriptsize $0101$}
        (112.5:2) node (9){\scriptsize $0010$}
        (135.0:2) node (10){\scriptsize $0011$}
        (157.5:2) node (11){\scriptsize $1001$}
        (180.0:2) node (12){\scriptsize $1110$}
        (202.5:2) node (13){\scriptsize $1011$}
        (225.0:2) node (14){\scriptsize $0100$}
        (247.5:2) node (15){\scriptsize $1000$}
        (270.0:2) node (1){\scriptsize $0111$}
        (292.5:2) node (3){\scriptsize $1100$}
        (315.0:2) node (5){\scriptsize $0000$}
        (337.5:2) node (6){\scriptsize $0110$};
      \begin{scope}[-]
        \draw (0) to (1);
        \draw (0) to (3);
        \draw (0) to (5);
        \draw (0) to (6);
        \draw (2) to (4);
        \draw (2) to (14);
        \draw (2) to (15);
        \draw (2) to (1);
        \draw (4) to (7);
        \draw (4) to (14);
        \draw (4) to (5);
        \draw (7) to (9);
        \draw (7) to (13);
        \draw (7) to (3);
        \draw (8) to (11);
        \draw (8) to (13);
        \draw (8) to (1);
        \draw (8) to (5);
        \draw (9) to (10);
        \draw (9) to (12);
        \draw (9) to (5);
        \draw (10) to (11);
        \draw (10) to (12);
        \draw (10) to (6);
        \draw (11) to (13);
        \draw (11) to (14);
        \draw (12) to (15);
        \draw (12) to (1);
        \draw (13) to (15);
        \draw (14) to (6);
        \draw (15) to (3);
        \draw (3) to (6);
      \end{scope}
    \end{tikzpicture}
    \caption{A $4$-regular closed XOR-magic graph of order $16$}
    \label{4regorder16}
\end{figure}
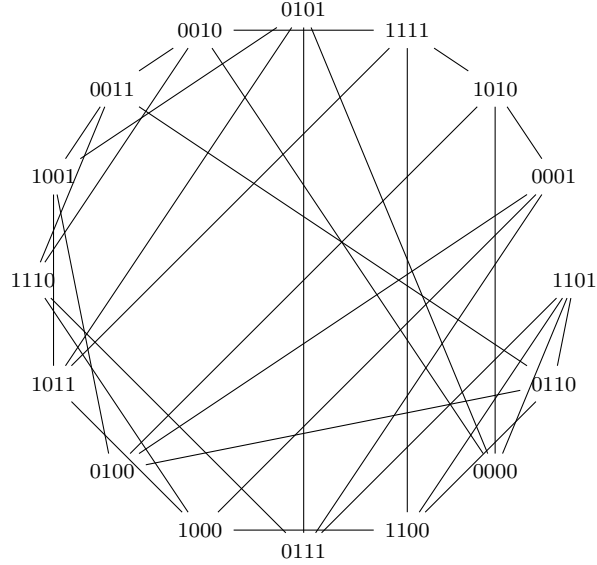

In particular, a $4$-regular closed XOR-magic graph of order $16$ is depicted in Figure \ref{4regorder16}.

\section{Conclusions}

An interesting observation comes from another fact proven by Batal. He proved the following fact.
\begin{proposition}[\cite{Batal}]\label{propbatalnieistnieja}
    There does not exist any even $k$-regular graph of order $2^n$ admitting a closed XOR-magic labeling for $k \in \{0,2,2^n-4,2^n-2\}$. Equivalently, there does not exist any odd $k$-regular graph of order $2^n$ admitting an open XOR-magic labeling for $k \in \{1,3,2^n-3, 2^n-1\}$.
\end{proposition}

Firstly, notice that the direct consequence of the Proposition \ref{propbatalnieistnieja} is that there does not exist any odd (even) regular open (closed) XOR-magic graph of order $8$. So, our constructed graphs of order $16$ are the minimal possible graphs. However, we proved that for order $2^n$ (where $n = 4$), there exist odd regular open XOR-magic graphs for each possible odd degree $d$ such that $d \notin \{1,3,2^n-3, 2^n-1\}$ and there exist even regular closed XOR-magic graphs for each possible even degree $d$ such that $d \notin \{0,2,2^n-2,2^n\}$.

Moreover, observe that by Theorems~\ref{openXORoddregularfig} and \ref{closedxomagicevenregular} and Corollary~\ref{cartesian} we obtain that there exists an open odd XOR-magic graph of power $8$ for any odd degree $d \in \{9,11,\ldots,21\}$ as well as an even closed XOR-magic graph of power $8$ for any odd degree $d \in \{8,10,\ldots,24,34,44,48,54,62,76,80,98,120\}$. 

\vspace{.3cm}
Recall, if $\delta(G)\geq (m-1)/2$ for a graph $G$ of order $m$, then the graph $G$ is connected. Thus, by Observation~\ref{complement} there exists odd $d$-regular open XOR-magic graph of power $8$ for any odd $d \in \{135,157,175,179,193,201,207,211,\allowbreak 221,231,\ldots,247\}$ as well as an even closed XOR-magic graph of power $8$ for any even degree $d \in \{234,236,\ldots,246\}$.
 In connection with that, we pose the following conjecture.

\begin{conjecture}
    For each degree $d \in \{5,7,\ldots 2^n-7, 2^n-5\}$, there exists $d$-regular open XOR-magic graph of order $2^n$. For each degree $d \in \{4,6,\ldots 2^n-6, 2^n-4\}$, there exists $d$-regular closed XOR-magic graph of order $2^n$.
\end{conjecture}

\end{document}